\documentclass[12pt]{amsart}
\usepackage[latin1]{inputenc}
\usepackage{amscd}
\usepackage{latexsym}
\usepackage{pstcol,pst-plot,pst-3d}
\usepackage{mathptmx}
\usepackage{multicol}

\psset{unit=0.7cm,linewidth=0.8pt,arrowsize=2.5pt 4}

\newpsstyle{fatline}{linewidth=1.5pt}
\newpsstyle{fyp}{fillstyle=solid,fillcolor=verylight}
\definecolor{verylight}{gray}{0.97}
\definecolor{light}{gray}{0.9}
\definecolor{medium}{gray}{0.85}

\def\xb{{\bold x}}
\def\yb{{\bold y}}

\def\opn#1#2{\def#1{\operatorname{#2}}} 
\opn\chara{char} \opn\length{\ell} \opn\pd{pd} \opn\rk{rk}
\opn\projdim{proj\,dim} \opn\injdim{inj\,dim} \opn\rank{rank}
\opn\depth{depth} \opn\grade{grade} \opn\height{height}
\opn\embdim{emb\,dim} \opn\codim{codim}

\opn\Tr{Tr} \opn\bigrank{big\,rank}
\opn\superheight{superheight}\opn\lcm{lcm}
\opn\trdeg{tr\,deg}%
\opn\reg{reg} \opn\lreg{lreg} \opn\skel{skel}
\opn\multideg{multideg}
\opn\div{div} \opn\Div{Div} \opn\cl{cl} \opn\Cl{Cl}
\opn\Spec{Spec} \opn\Supp{Supp} \opn\supp{supp} \opn\Sing{Sing}
\opn\Ass{Ass}
\opn\Ann{Ann} \opn\Rad{Rad} \opn\Soc{Soc}
\opn\Ker{Ker} \opn\Coker{Coker} \opn\Im{Im} \opn\Hom{Hom}
\opn\Tor{Tor} \opn\Ext{Ext} \opn\End{End} \opn\Aut{Aut}
\opn\id{id}

\opn\nat{nat}
\opn\pff{pf}
\opn\Pf{Pf} \opn\GL{GL} \opn\SL{SL} \opn\mod{mod} \opn\ord{ord}
\opn\aff{aff} \opn\con{conv} \opn\relint{relint} \opn\st{st}
\opn\lk{lk} \opn\cn{cn} \opn\core{core} \opn\vol{vol}
\opn\link{link} \opn\star{star} \opn\skel{skel}
\opn\gr{gr}

\def\pot#1#2{#1[\kern-0.28ex[#2]\kern-0.28ex]}

\opn\dirlim{\underrightarrow{\lim}}
\opn\inivlim{\underleftarrow{\lim}}
\let\union=\cup

\let\to=\rightarrow

\def\Implies{\ifmmode\Longrightarrow \else
     \unskip${}\Longrightarrow{}$\ignorespaces\fi}
\def\implies{\ifmmode\Rightarrow \else
     \unskip${}\Rightarrow{}$\ignorespaces\fi}
\def\iff{\ifmmode\Longleftrightarrow \else
     \unskip${}\Longleftrightarrow{}$\ignorespaces\fi}

\let\:=\colon
\newtheorem{Theorem}{Theorem}[section]
\newtheorem{Lemma}[Theorem]{Lemma}
\newtheorem{Corollary}[Theorem]{Corollary}
\newtheorem{Proposition}[Theorem]{Proposition}

\newtheorem{Example}[Theorem]{Example}

\newtheorem{Question}[Theorem]{Question}

\let\epsilon\varepsilon
\let\phi=\varphi
\let\kappa=\varkappa
\def\qed{\ifhmode\textqed\fi
   \ifmmode\ifinner\quad\qedsymbol\else\dispqed\fi\fi}
\def\textqed{\unskip\nobreak\penalty50
    \hskip2em\hbox{}\nobreak\hfil\qedsymbol
    \parfillskip=0pt \finalhyphendemerits=0}
\def\dispqed{\rlap{\qquad\qedsymbol}}

\def\BB{{\mathcal B}}
\def\FF{{\mathcal F}}

\def\LL{{\mathcal L}}
\def\II{{\mathcal I}}
\def\JJ{{\mathcal J}}

\opn\inii{in} \opn\inim{inm} \opn\rate{rate}

\begin{document}
\title{Level rings arising from meet-distributive meet-semilattices}
\author{J\"urgen Herzog and  Takayuki Hibi}

\address{J\"urgen Herzog, Fachbereich Mathematik und
Informatik, Universit\"at Duisburg-Essen, Campus Essen,
 45117 Essen, Germany}
\email{juergen.herzog@uni-essen.de}
\address{Takayuki Hibi, Department of Pure and Applied Mathematics,
 Graduate School of Information Science and Technology,
Osaka University, Toyonaka, Osaka 560-0043, Japan}
\email{hibi@math.sci.osaka-u.ac.jp}  \maketitle
\begin{abstract}
The Alexander dual of an arbitrary meet-semilattice is described
explicitly.  Meet-distributive meet-semilattices whose Alexander
dual is level are characterized.
\end{abstract}

\section*{Introduction}
In the present paper we continue our discussion in \cite{HerHi}
and \cite{HHZHHZ}, and describe explicitly the generators of the
Alexander dual of an arbitrary meet-semilattice (Theorem
\ref{AlexanderDual}). In case of meet-distributive
meet-semilattices, a combinatorial formula (Proposition
\ref{hvector}) to compute the $h$-vector of the Alexander dual is
given.

It is known \cite{HHZHHZ} that the Alexander dual $\Gamma_\LL$ of
a meet-semilattice $\LL$ is Cohen--Macaulay if and only if $\LL$
is meet-distributive.  Thus it seems of interest to characterize
the meet-distributive meet-semilattices $\LL$ for which
$\Gamma_\LL$ is a level complex \cite{StanleySecondEd}.  Our main
theorem (Theorem \ref{maintheorem}) says that the Alexander dual
$\Gamma_\LL$ is level if and only if a certain simplicial complex
coming from $\LL$ is pure.  In particular, in case that $\LL$ is a
finite distributive lattice, $\Gamma_\LL$ is level if and only if
the simplicial complex consisting of all antichains of the poset
of all join-irreducible elements of $\LL$ is pure.

\section{The $h$-vector of a finite meet-distributive
meet-semilattice} First of all, we prepare notation and
terminologies on finite lattices and finite posets (partially
ordered sets). In a finite poset $P$ we say that $\alpha \in P$
{\em covers} $\beta \in P$ (or $\beta$ is a {\em lower neighbor}
of $\alpha$) if $\beta < \alpha$ and $\beta < \gamma < \alpha$ for
no $\gamma \in P$. Let $N(\alpha)$ denote the set of lower
neighbors of $\alpha \in P$. A {\em poset ideal} of $P$ is a
subset $\II$ of $P$ such that $\alpha \in \II$ and $\beta \in P$
together with $\beta \leq \alpha$ imply $\beta \in \II$.

Let $\LL$ be a finite meet-semilattice \cite[p.\ 103]{StanleyEC}
%
%
%
%
and ${\hat 0}$ its unique minimal element.
Since $\LL$ is a meet-semilattice,
it follows from
\cite[Proposition 3.3.1]{StanleyEC}
that $\LL$ is a lattice if and only if
$\LL$ possesses a unique maximal element ${\hat 1}$.
In other words, if $\LL$ is a meet-semilattice
and is not a lattice, then $\LL \cup \{ {\hat 1} \}$
with a new element ${\hat 1}$
such that $\alpha < {\hat 1}$
for all $\alpha \in \LL$ becomes a lattice.
Thus, in a finite meet-semilattice $\LL$,
each element of $\LL$ is the join of elements
of $\LL$.
A {\em join-irreducible element} of $\LL$
is an element $\alpha \in \LL$
such that one cannot write
$\alpha  = \beta \vee \gamma$
with $\beta < \alpha$ and $\gamma < \alpha$.
In other words, a join-irreducible element
of $\LL$
is an element $\alpha \in \LL$
which covers exactly one element of $\LL$.

Let $\LL$ be a finite meet-semilattice
and $P \subset \LL$ the set of join-irreducible
elements of $\LL$.
We will associate each element $\alpha \in \LL$
with the subset
\begin{eqnarray}
\label{canonicalembedding}
\ell(\alpha) = \{ p \in P \, : \, p \leq \alpha \}.
\end{eqnarray}
Thus $\ell(\alpha)$ is a poset ideal of $P$,
and $\alpha \in \ell(\alpha)$ if and only if
$\alpha$ is join-irreducible.
Moreover, for $\alpha$ and $\beta$
belonging to $\LL$, one has
$\ell(\alpha) = \ell(\beta)$
if and only if $\alpha = \beta$.

\begin{Lemma}
\label{canonicalembedding}
One has
$\ell(\alpha \wedge \beta)
= \ell(\alpha) \cap \ell(\beta)$
for all $\alpha, \beta \in \LL$.
\end{Lemma}

\begin{proof}
Let
$\gamma = \alpha \wedge \beta$.
Then
$\ell(\gamma) \subset
\ell(\alpha) \cap \ell(\beta)$.
Since $\LL \cup \{ {\hat 1} \}$
with a new element ${\hat 1}$
is a lattice,
if
$\ell(\gamma) \neq
\ell(\alpha) \cap \ell(\beta)$
and if $p \in
(\ell(\alpha) \cap \ell(\beta))
\setminus \ell(\gamma)$,
then
$\delta = \gamma \vee p \in \LL$
with $\gamma < \delta \leq \alpha$
and $\delta \leq \beta$.
This contradicts
$\gamma = \alpha \wedge \beta$.
\end{proof}

Let $K$ be a field and
$K[\xb, \yb] = K[\{ x_p, y_p \}_{p \in P}]$
denote the polynomial
ring in $2|P|$ variables over $K$
with each $\deg x_p = \deg y_p = 1$.
We associate each element $\alpha \in \LL$
with the squarefree monomial
\[
u_\alpha = (\prod_{p \in \ell(\alpha)} x_p)
(\prod_{p \in P \setminus \ell(\alpha)} y_p) \in
K[\xb, \yb]
\]
and set
\[
H_{\LL} = (u_\alpha)_{\alpha \in \LL}
\subset K[\xb, \yb].
\]
Since the ideal $H_{\LL}$ is squarefree,
there is a simplicial complex $\Sigma_\LL$
on the vertex set
$\{ x_p, y_p \}_{p\in P}$ whose
Stanley--Reisner ideal $I_{\Sigma_\LL}$ coincides
with $H_{\LL}$.

Let $\Gamma_\LL$ denote the {\em Alexander dual}
(\cite{EagonReiner}, \cite{HHZDirac}) of $\Sigma_\LL$
%
%
%
%
and call $\Gamma_\LL$
the Alexander dual of
$\LL$. We write ${\mathcal F}(\Gamma_\LL)$ for the set of facets
(maximal faces) of $\Gamma_\LL$. One has
\begin{eqnarray}
\label{facets} {\mathcal F}(\Gamma_\LL) = \{F_\alpha \, : \,
\alpha \in \LL\}
\end{eqnarray}
where
\[
F_\alpha
= \{ x_q \, : \, q \in P \setminus \ell(\alpha) \}
\union \{ y_q \: q \in  \ell(\alpha)\}.
\]
In particular
$\Gamma_\LL$ is a pure simplicial complex of dimension
$|P| - 1$.

A finite meet-semilattice $\LL$ is called
{\em meet-distributive}
\cite[p. 156]{StanleyEC}
if each interval $[\alpha, \beta]
= \{ \gamma \in \LL \:
\alpha \leq \gamma \leq \beta \}$
of $\LL$ such that $\alpha$
is the meet of the lower neighbors
of $\beta$ in $[\alpha, \beta]$ is boolean.
For example, every poset ideal of a finite
distributive lattice is a meet-distributive
meet-semilattice.

Let $\LL$ be an arbitrary finite meet-distributive
meet-semilattice and, as before, $P \subset \LL$
the set of join-irreducible elements of $\LL$.
The {\em distributive closure} of $\LL$ is
the finite distributive lattice $\JJ(P)$
consisting of all poset ideals of $P$
ordered by inclusion.

Recall that Birkhoff's fundamental structure theorem on finite
distributive lattices \cite[Theorem 3.4.1]{StanleyEC} guarantees
that every finite distributive lattice is of the form $\JJ(P)$ for
a unique finite poset $P$. In fact, if $P$ is the set of
join-irreducible element of a finite distributive lattice $\LL$,
then $\LL = \JJ(P)$.

It is not difficult to see that the map
   $\ell \: \LL \to \JJ(P)$ defined by
   (\ref{canonicalembedding}) is an embedding
   of meet-semilattices if and only if $\LL$ is
   meet-distributive.
   Consult \cite{Edelman} for further information
   about meet-distributive lattice.

\begin{Proposition}
\label{hvector}
Let $\LL$ be a finite meet-distributive meet-semilattice
and $\Gamma_\LL$ its Alexander dual.
Let $h(\Gamma_\LL) = (h_0, h_1, \ldots)$
be its $h$-vector.  Then, for all $i$, one has
\[
h_i = |\{\alpha\in\LL \: |N(\alpha)| = i \}|.
\]
\end{Proposition}

\begin{proof}
Let $\alpha \in \LL$ with $|N(\alpha)|=i$ and $\ell(\alpha) = \{
q_1, \ldots, q_\delta \}$. Let $N(\alpha) = \{ r_1, \ldots, r_i
\}$ with each $\ell(r_j) = \ell(\alpha) \setminus \{ q_{\delta - j
+ 1} \}$. Let $r = r_1 \wedge \cdots \wedge r_i$. Thus $\ell(r) =
\cap_{j=1}^{i} \ell(r_j) = \{ q_1, \ldots, q_{\delta - i} \}$ and
the interval $[r, \alpha]$ in $\LL$ is the boolean lattice of rank
$i$. Since a subset $A \subset \ell(\alpha)$ is contained in none
of the sets $\ell(r_1), \ldots, \ell(r_i)$ if and only if $A$
contains $\{ q_\delta , q_{\delta - 1}, \ldots, q_{\delta - i + 1}
\}$, it follows that the number of subsets $A \subset
\ell(\alpha)$ with $|A| = k$ such that $A \subset \ell(q)$ for no
$q \in \LL$ with $q < \alpha$ is ${ \delta - i \choose k - i }$.
In other words, the number of those faces $F \subset F_\alpha$ of
$\Gamma$ with $|F| = j + 1$ such that $F \subset F_q$ for no $q
\in \LL$ with $q < \alpha$ is $\sum_{k=i}^{\delta} { |P| - \delta
\choose j - k + 1}{ \delta - i \choose k - i }$, which is equal to
${ |P| - i \choose j - i + 1 } = { |P| - i \choose |P| - j - 1 }$.
Thus the number of faces $F$ of $\Gamma_\LL$ with $|F| = j + 1$ is
\[
f_j(\Gamma_\LL) = \sum_{i=0}^{j+1} { |P| - i \choose |P| - j - 1 }
|\{\alpha\in\LL\: |N(\alpha)|=i\}|.
\]
On the other hand, in general, one has
\[
f_j(\Gamma_\LL) = \sum_{i=0}^{j+1} { |P| - i \choose |P| - j - 1 }
h_i(\Gamma_\LL).
\]
Hence $h_i(\Gamma_\LL) = |\{\alpha\in\LL\: |N(\alpha)|=i\}|$, as
desired.
\end{proof}

\begin{Corollary}
\label{easy} Let $\LL$ be a finite meet-distributive
meet-semilattice, $P$ the set of join-irreducible elements of
$\LL$ and $\Gamma_\LL$ its Alexander dual. Let $n = |P|$ and
$(h_0, h_1, \ldots, h_n)$ the $h$-vector of $\Gamma_\LL$. Then
$h_1 = n$, and the $a$-invariant of $\Gamma_\LL$ (which is the
nonpositive integer $n - \max\{ i \: h_i \neq 0 \}$) is equal to
$\max\{ |N(\alpha)| \: \alpha \in\LL \} - n$.
\end{Corollary}

\begin{Example}
\label{polarization} {\em Let $\BB_{[n]}$ denote the boolean
lattice of rank $n$ and $\LL$ a poset ideal of $\BB_{[n]}$ which
contains all join-irreducible elements (i.e., $\{ 1 \}, \ldots, \{
n \}$) of $\BB_{[n]}$.   Then the meet-distributive
meet-semilattice $\LL$ is a simplicial complex on $[n]$ and the
$h$-vector of $\Gamma_\LL$ coincides with the $f$-vector of $\LL$.

(a)
By using (\ref{facets})
the Stanley--Reisner ideal $I_{\Gamma_\LL}$ of $\Gamma_\LL$ is
generated by those squarefree monomials $\prod_{q \in \ell(\beta)}
y_q$ such that $\beta \in \BB_{[n]}$ is a minimal nonface of $\LL$
and by the quadratic monomials $x_{\{ i \}} y_{\{ i \}}$ for all
$i \in [n]$.

(b) Let $T = K[y_{\{ 1 \}}, \ldots, y_{\{ n \}}]$ and $J \subset
T$ the ideal generated by those squarefree monomials $\prod_{q \in
\ell(\beta)} y_q$ such that $\beta \in \BB_{[n]}$ is a minimal
nonface of $\LL$ and by $y_{\{ i \}}^2$ for all $i \in [n]$. The
quotient ring $T / J$ is of $0$-dimensional and its $h$-vector is
$(f_{-1}, f_0, \ldots)$, the $f$-vector of $\LL$ with $f_{-1} =
1$. It turns out that $I_{\Gamma_\LL}$ is the polarization
\cite[Lemma 4.3.2]{BrunsHer} of the ideal $J$.
%
%
%
%
Since $T / J$ is Cohen--Macaulay,
it follows immediately that $\Gamma_\LL$
is Cohen--Macaulay.
This fact is a special case of
\cite[Corollary 1.6]{HHZHHZ}.

(c)
Since $T / J$
is a level ring \cite[p. 91]{StanleySecondEd}
%
%
%
%
if and only if the simplicial complex $\LL$
is pure, it follows that the Alexander dual
$\Gamma_\LL$ of $\LL$
is a level complex if and only if
the simplicial complex $\LL$ is pure.

(d) Let $\Delta$ be a simplicial complex on the vertex set $V = \{
y_{\{ 1 \}}, \ldots, y_{\{ n \}} \}$, and let $W = \{ x_{\{ 1 \}},
\ldots, x_{\{ n \}}\}$. We write $\Delta^\sharp$ for the
simplicial complex on the vertex set $V \cup W$ whose facets are
those of $\Delta$ together with all edges $\{ x_{\{ i \}}, y_{\{ i
\}} \}$ for $i = 1, \ldots, n$. By the observation (a) for a
simplicial complex $\LL$ $( \subset \BB_{[n]} )$ on $[n]$ one has
a simplicial complex $\Delta$ on $V$ such that the facet ideal of
$\Delta^\sharp$, i.e., the ideal generated by all monomials
corresponding to the facets, coincides with the Stanley--Reisner
ideal $I_{\Gamma_\LL}$ of $\Gamma_\LL$. Conversely, given a
simplicial complex $\Delta$ on $V$, there is a simplicial complex
$\LL$ $( \subset \BB_{[n]} )$ on $[n]$ such that the facet ideal
of $\Delta^\sharp$ coincides with $I_{\Gamma_\LL}$. Since
$\Gamma_\LL$ is always Cohen--Macaulay, the facet ideal of
$\Delta^\sharp$ is Cohen--Macaulay. This argument is a direct and
easy proof of \cite[Corollary 4.4]{HHZHHZ}.

(e) By using (b) and (c), it follows that
every $f$-vector of a pure simplicial complex
is the $h$-vector of a level complex.
}
\end{Example}

It would, of course, be of interest to generalize
the fact (c) of Example \ref{polarization}
to arbitrary meet-distributive meet-semilattices
$\LL$.


\section{Alexander duality of meet-distributive
meet-semilattices}
A nice description of the Alexander dual of a finite
distributive lattice is obtained in
\cite[Lemma 3.1]{HerHi}.
%
%
%
%
On the other hand, the Alexander dual of
a meet-distributive meet-semilattice of a special kind,
namely, a poset ideal of a finite distributive lattice
is described in \cite[Theorem 4.2]{HHZHHZ}.
An explicit description of the Alexander dual
of an arbitrary finite meet-semilattice
will be obtained
in Theorem \ref{AlexanderDual} below.

If, in general, $P$ is a finite poset
and $B \subset P$, then we write
$\langle B \rangle$ for the poset ideal of $P$
{\em generated by $B$}, i.e.,
$p \in P$ belongs to $\langle B \rangle$
if and only if $p \leq q$ for some $q \in B$.

\begin{Theorem}
\label{AlexanderDual} Let $\LL$ be an arbitrary finite
meet-semilattice and $P$ the set of join-irreducible elements of
$\LL$. Then the Stanley--Reisner ideal $I_{\Gamma_\LL}$ of the
Alexander dual $\Gamma_\LL$ of $\LL$ is generated by the following
squarefree monomials:
\begin{enumerate}
\item[(i)] $x_p y_q$, where $p, q \in P$ with $p < q$; \item[(ii)]
$\prod_{q \in B} y_q$, where $B$ is an antichain of $P$ with
$\langle B \rangle \not \subset \ell(\alpha)$ for all $\alpha \in
\LL$; \item[(iii)] $x_p \prod_{q \in B} y_q$, where
$B$ is an antichain of $P$
with $\ell(\beta) \neq \langle B \rangle$ for all $\beta \in \LL$,
but with $\langle B \rangle \subset \ell(\alpha)$ for some $\alpha
\in \LL$ and where
$p \in \ell(\wedge_{\langle B \rangle
\subset \ell(\alpha)}
\alpha) \setminus \langle B \rangle$.
\end{enumerate}
\end{Theorem}

\begin{proof}
Let $A \subset P$ and $B \subset P$ with
$A \cap B = \emptyset$.
We write $\xb_A \yb_B$
for the squarefree monomial
$\prod_{p \in A} x_p \prod_{q \in B} y_q$
of $K[\xb, \yb]$.
By the definition of the Stanley--Reisner ideal
$I_{\Gamma_\LL}$ of $\Gamma_\LL$,
it follows that
$\xb_A \yb_B$ belongs to
$I_{\Gamma_\LL}$ if and only if
there is no facet $\FF_\alpha$
of $\Gamma_\LL$
with
$\{ x_p \, : \, p \in A\}
\cup \{ y_q \, : \, q \in B \}
\subset \FF_\alpha$.
Thus by using (\ref{facets})
one has
$\xb_A \yb_B \in I_{\LL_\Gamma}$
if and only if
there is no $\alpha \in \LL$
such that
$A \subset P \setminus \ell(\alpha)$
and
$B \subset \ell(\alpha)$.
In other words, one has
$\xb_A \yb_B \in I_{\LL_\Gamma}$
if and only if
the following condition
$(*)$ is satisfied:
\[
\text{
$(*)$
\, \, \, \, \,
each $\alpha \in \LL$
with
$B \subset \ell(\alpha)$
satisfies
$A \cap \ell(\alpha) \neq \emptyset$.
\, \, \, \, \, \, \, \, \, \,
}
\]
We say that a pair
$(A, B)$, where
$A \subset P$ and $B \subset P$ with
$A \cap B = \emptyset$,
is an {\em independent pair} of $\LL$
if the condition $(*)$ is satisfied.
Thus
$I_{\LL_\Gamma}$
is generated by all monomials $x_py_p$
with $p \in P$ together with
those monomials $\xb_A \yb_B$
such that $(A, B)$
is an independent pair of $\LL$.


Let $M(B)$ denote the set of maximal elements of
$B$.  Thus one has
$\langle B \rangle = \langle M(B) \rangle$.
Hence $(A, B)$ is independent
if and only if
$(A, M(B))$ is independent.
Since $\xb_A \yb_{M(B)}$ divides $\xb_A \yb_B$
and since $M(B)$ is an antichain of $P$,
it follows that
$I_{\LL_\Gamma}$
is generated by all monomials $x_py_p$
with $p \in P$ together with
those monomials $\xb_A \yb_B$
such that $(A, B)$
is an independent pair of $\LL$
and $B$ is an antichain of $P$.


Let $p$ and $q$ belong to $P$.
Since $\ell(q) = \langle \{ q \} \rangle \in \LL$
for all $q \in \LL$,
the pair $(\{ p \}, \{ q \})$ with $p \neq q$
is an independent pair of $\LL$
if and only if
$p < q$.
Let $\ell(\beta) = \langle B \rangle$
for some $\beta \in \LL$.
Then a pair
$(A, B)$ is independent if and only if
$A \cap \langle B \rangle \neq \emptyset$.
On the other hand,
$A \cap \langle B \rangle \neq \emptyset$
if and only if
there are $p \in A$ and $q \in B$
with $p < q$.

Consequently,
$I_{\LL_\Gamma}$
is generated by all monomials $x_py_q$,
where $p, q \in P$ with $p < q$
together with
those monomials $\xb_A \yb_B$, where
$(A, B)$
is an independent pair of $\LL$
such that
$B$ is an antichain of $P$
with $\ell(\beta) = \langle B \rangle$
for no $\beta \in \LL$
and with $A \cap \langle B \rangle = \emptyset$.

Now, let $B$ be an antichain of $P$
with $\ell(\beta) = \langle B \rangle$
for no $\beta \in \LL$ and
$A \subset P$ with
$A \cap \langle B \rangle = \emptyset$.

(a)
First, if
$\langle B \rangle \subset \ell(\alpha)$
for no $\alpha \in \LL$, then
$(A, B)$ is independent
for all $A \subset P$ with $A \cap B = \emptyset$.
Thus in particular $(\emptyset, B)$
is an independent pair of $\LL$.

(b)
Second, if
$\langle B \rangle \subset \ell(\alpha)$
for some $\alpha \in \LL$, then
$(A, B)$ is independent if and only if
$A \cap (\cap_{\langle B \rangle \subset \ell(\alpha)}
\ell(\alpha)) \neq \emptyset$.
Since $\cap_{\langle B \rangle \subset \alpha}
\ell(\alpha)
= \ell(\wedge_{\langle B \rangle \subset \alpha}
\alpha)$,
it follows that
$(A, B)$ is independent if and only if
there is $p \in A$ with
$p \in \ell(\wedge_{\langle B \rangle \subset \alpha}
\alpha) \setminus \langle B \rangle$.
\end{proof}

\section{Level rings arising from Meet-distributive
meet-semilattices}
Recall that the Alexander dual $\Gamma_\LL$
of a finite meet-semilattice $\Gamma$ is
Cohen--Macaulay if and only if
$\LL$ is meet-distributive
(\cite[Corollary 1.6]{HHZHHZ}).
The problem when the Alexander dual of
a meet-distributive meet-semilattice is a level ring
is now studied.

Let $\LL$ be a finite meet-distributive meet-semilattice
and $P$ the set of join-irreducible elements
of $\LL$.  Let, as before,
$K[\xb, \yb] = K[\{ x_p, y_p \}_{p \in P}]$
denote the polynomial
ring in $2|P|$ variables over a field $K$
with each $\deg x_p = \deg y_p = 1$.

For each $\alpha \in \LL$ we write $\alpha' \in\LL$ for
the meet of all $\beta \in N(\alpha)$, where
$N(\alpha)$ is the set of lower neighbors of $\alpha$.
Since $\LL$ is a meet-distributive meet-semilattice,
it follows that the interval
\[
\BB_\alpha =
[\alpha', \alpha]
= \{ \gamma \in \LL \:
\alpha' \leq \gamma \leq \alpha \}
\]
of $\LL$ is a boolean lattice.
Let $S(\alpha) \subset P$ denote the antichain
\[
S(\alpha) = \ell(\alpha) \setminus \ell(\alpha').
\]
Each element belonging to $S(\alpha)$ is
a maximal element of $\ell(\alpha)$, and
$t \in \ell(\alpha)$ belongs to $S(\alpha)$
if and only if
$\ell(\beta) = \ell(\alpha) \setminus \{ t \}$
for some $\beta \in N(\alpha)$.

\begin{Lemma}
\label{differentfrom}
If $\alpha$ and $\beta$ belong to $\LL$
with $\alpha \neq \beta$,
then $S(\alpha) \neq S(\beta)$.
\end{Lemma}

\begin{proof}
Let $\gamma = \alpha \wedge \beta$.
If $S(\alpha) = S(\beta)$, then
$S(\alpha) \subset S(\gamma)$.
In fact, for each $t \in S(\alpha) = S(\beta)$,
there are $\alpha_0 \in N(\alpha)$ and
$\beta_0 \in N(\beta)$ with
$\ell(\alpha_0) = \ell(\alpha) \setminus \{ t \}$
and
$\ell(\beta_0) = \ell(\beta) \setminus \{ t \}$.
By using Lemma \ref{canonicalembedding},
one has
$\ell(\alpha_0 \wedge \beta_0)
= \ell(\gamma) \setminus \{ t \}$.
Hence $t \in S(\gamma)$.
Now, since $\gamma < \alpha$, one has
$\delta \in N(\alpha)$ with
$\gamma \leq \delta < \alpha$.
Since $\ell(\delta) = \ell(\alpha)
\setminus \{ t \}$ for some $t \in S(\alpha)$,
it follows that $t \not\in \ell(\gamma)$.
This contradict $S(\alpha) \subset S(\gamma)$.
Hence $S(\alpha) \neq S(\beta)$, as desired.
\end{proof}

Recall from the proof of Theorem \ref{AlexanderDual}
that a pair $(A, B)$, where
$A \subset P$ and $B \subset P$ with
$A \cap B = \emptyset$,
is said to be an independent pair of $\LL$
if each $\alpha \in \LL$ with
$B \subset \ell(\alpha)$
satisfies
$A \cap \ell(\alpha) \neq \emptyset$.

\begin{Lemma}
\label{monomialbasis}
Let $\alpha \in \LL$
and $T \subset S(\alpha)$.
Then the pair $(\emptyset, T)$
cannot be independent.
Moreover,
for $p \in S(\alpha) \setminus T$,
the pair $(\{ p \}, T)$
cannot be independent.
\end{Lemma}

\begin{proof}
Since $T \subset \ell(\alpha)$,
the pair $(\emptyset, T)$
cannot be an independent pair of $\LL$.
On the other hand,
since $p \in S(\alpha)$,
one has $\beta \in N(\alpha)$
with
$\ell(\beta) = \ell(\alpha) \setminus \{ p \}$.
Since $T \subset \ell(\beta)$
and since $\{ p \} \cap \ell(\beta) = \emptyset$,
it follows that $(\{ p \}, T)$
cannot be an independent pair of $\LL$,
as desired.
\end{proof}

Let $I_{\Gamma_\LL}$ denote the Stanley--Reisner ideal
of
$\LL$ and
$K[\Gamma_\LL] = K[\xb, \yb] / I_{\Gamma_\LL}$
the Stanley--Reisner ring of $\Gamma_\LL$.
Since the dimension of $\Gamma_\LL$ is
$|P| - 1$ and the Krull dimension of $K[\Gamma_\LL]$
coincides with $|P|$,
it follows easily that
$\{ x_p - y_p \, : \, p \in P \}$
is a linear system of parameters of $K[\Gamma_\LL]$.
Since $K[\Gamma_\LL]$ is Cohen--Macaulay, by using
Proposition \ref{hvector}, the Hilbert series of
the quotient ring
\[
K[\Gamma_\LL] / (x_p - y_p \, : \, p \in P)
\]
is $h_0 + h_1 \lambda + h_2 \lambda^2 + \cdots$, where $(h_0, h_1,
h_2, \ldots)$ is the $h$-vector of $\Gamma_\LL$.

Let $J_{\Gamma_\LL}$ be the monomial ideal of $K[\xb] = K[\{ x_p
\}_{p \in P}]$ generated by those monomials
\begin{enumerate}
\item[(i)]
$x_p x_q$, where
$p, q \in P$ with $p < q$;
\item[(ii)]
$\prod_{q \in B} x_q$, where
$B$ is an antichain of $P$
with $\langle B \rangle \subset \ell(\alpha)$
for no $\alpha \in \LL$;
\item[(iii)]
$x_p \prod_{q \in B} x_q$,
where
$B$ is an antichain of $P$
with $\ell(\beta) = \langle B \rangle$
for no $\beta \in \LL$, but
with
$\langle B \rangle \subset \ell(\alpha)$
for some $\alpha \in \LL$
and where
$p \in \ell(\wedge_{\langle B \rangle
\subset \ell(\alpha)}
\alpha) \setminus \langle B \rangle$.
\end{enumerate}
By virtue of Theorem \ref{AlexanderDual}
it follows that
\[
K[\xb] / J_{\Gamma_\LL}
=
K[\Gamma_\LL] / (x_p - y_p \, : \, p \in P).
\]
We associate each $\alpha \in \LL$ with
the monomial
\[
u_\alpha = \prod_{p \in S(\alpha)} x_p
\]
of degree $|N(\alpha)|$.

\begin{Theorem}
\label{maintheorem}
Let $\LL$ be a finite meet-distributive meet-semilattice,
$P$ the set of join-irreducible elements of $\LL$,
and $K[\xb] = K[\{ x_p \}_{p \in P}]$
the polynomial ring in $|P|$ variables over a field $K$.
Then the set of monomials
$\{ u_\alpha \, ; \, \alpha \in \LL \}$
is a $K$-basis of the quotient ring
$K[\xb] / J_{\Gamma_\LL}$.
Thus in particular
$\{ S(\alpha) \, : \, \alpha \in \LL \}$
is a simplicial complex on the vertex set
$\{ x_p \, : \, p \in P \}$
whose $f$-vector coincides with
the $h$-vector of $\LL$.
\end{Theorem}

\begin{proof}
Lemma \ref{monomialbasis}
says that, for each $\alpha \in \LL$,
the monomial $u_\alpha$ does not belongs to
$J_{\Gamma_\LL}$.
Moreover, Lemma \ref{differentfrom}
guarantees that, for $\alpha \neq \beta$
belonging to $\LL$,
one has $u_\alpha \neq u_\beta$.
Hence, for each $i = 0, 1, 2, \ldots$,
the number of monomials $u_\alpha$
with $\alpha \in \LL$ of degree $i$
is equal to $h_i$.
Since the Hilbert series of
$K[\xb] / J_{\Gamma_\LL}$ is
$h_0 + h_1 \lambda + h_2 \lambda^2 + \cdots$,
it follows that
$\{ u_\alpha \, ; \, \alpha \in \LL \}$
is a $K$-basis of
$K[\xb] / J_{\Gamma_\LL}$,
as required.
\end{proof}

We now come to the combinatorial characterization
for the Alexander dual $\Gamma_\LL$ of
a finite meet-distributive meet-semilattice $\LL$
to be level.

\begin{Corollary}
\label{levelring}
The Alexander dual $\Gamma_\LL$ of
a finite meet-distributive meet-semilattice $\LL$
is a level complex if and only if
the simplicial complex
$\{ S(\alpha) \, : \, \alpha \in \LL \}$
is pure.
Thus in particular the
Alexander dual $\Gamma_\LL$ of
a finite distributive lattice $\LL = \JJ(P)$
is level if and only if
the simplicial complex consisting of all
antichains of $P$ is pure.
\end{Corollary}

\newpage
Consider the following example of a meet-distributive
meet-semilattice $\LL$
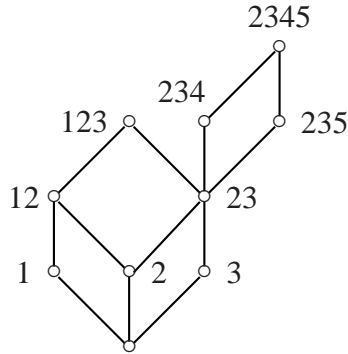
\begin{center}
 \psset{unit=2cm}
\begin{pspicture}(-2.5,-0.5)(5.5,2.3)
\psline(0.97,0.03)(0.53,0.47)
 \psline(1.47,0.47)(1.03,0.03)
 \psline(0.5,0.54)(0.5,0.97)
 \psline(0.53,1.03)(0.97,1.47)
 \psline(1.03,1.47)(1.47,1.03)
 \psline(1.5,0.97)(1.5,0.54)
 \psline(1,0.04)(1,0.47)
 \psline(0.53,0.97)(0.97,0.53)
 \psline(1.03,0.5)(1.47,0.97)
\psline(1.53,1.03)(1.97,1.47)

\psline(2,1.54)(2,1.97)

\psline(1.97,1.97)(1.53,1.53)
 \psline(1.5,1.04)(1.5,1.47)
\rput(2, 2){$\circ$}
 \rput(2,1.5){$\circ$}
 \rput(1,0){$\circ$}
 \rput(0.5,0.5){$\circ$}
 \rput(1,0.5){$\circ$}
 \rput(1.5,1.5){$\circ$}
 \rput(0.5,1){$\circ$}
 \rput(1.5,0.5){$\circ$}
 \rput(1.5,1){$\circ$}
 \rput(1,1.5){$\circ$}
 \rput(0.3,0.5){$1$}
 \rput(1.2,0.5){$2$}
 \rput(1.7,0.5){$3$}
 \rput(0.3,1){$12$}
 \rput(1.75,1){$23$}
 \rput(0.7,1.5){$123$}
 \rput(1.35,1.7){$234$}
 \rput(2.3,1.5){$235$}
 \rput(2,2.2){$2345$}
\end{pspicture}
\end{center}

\medskip
\noindent
with the following poset of join-irreducible elements
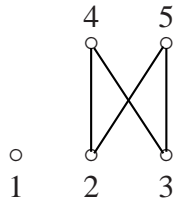
\begin{center}
 \psset{unit=2cm}
\begin{pspicture}(-2.5,-0.5)(5.5,1.2)
\rput(1,0){$\circ$} \rput(0.5,0){$\circ$} \rput(1.5,0){$\circ$}
 \rput(1,0.75){$\circ$} \rput(1.5,0.75){$\circ$}

\psline(1,0.03)(1,0.72) \psline(1.5,0.03)(1.5,0.72)
\psline(1.02,0.02)(1.48,0.73) \psline(1.02,0.72)(1.48,0.02)

\rput(1,-0.2){2}\rput(0.5,-0.2){1}\rput(1.5,-0.2){3}
\rput(1,0.93){4}\rput(1.5,0.93){5}
\end{pspicture}
\end{center}

By using Theorem \ref{AlexanderDual} the Stanley--Reisner ideal of
the Alexander dual of $\LL$ is generated by the following
monomials:
\begin{enumerate}
\item[(i)] $x_1y_1,\ldots, x_5y_5$, $x_2y_4$, $x_2y_5$, $x_3y_4$
and $x_3y_5$;

\item[(ii)] $y_1y_4$ and $y_1y_5$;

\item[(iii)] $x_2y_1y_3$.
\end{enumerate}

The $h$-vector of $\LL$ is $(1,5,4)$.  By using Corollary
\ref{levelring} the Alexander dual $\Gamma_\LL$ is level.

\medskip
\noindent In the following meet-distributive meet-semilattice
$\LL$ the facets of the simplicial complex $\{S(\alpha)\:
\alpha\in \LL\}$ are $\{1,2\}$, $\{1,3,4\}$ and $\{2,3,4\}$. Since
this simplicial complex is not pure it follows from Corollary
\ref{levelring} that $\Gamma_\LL$ is not level.
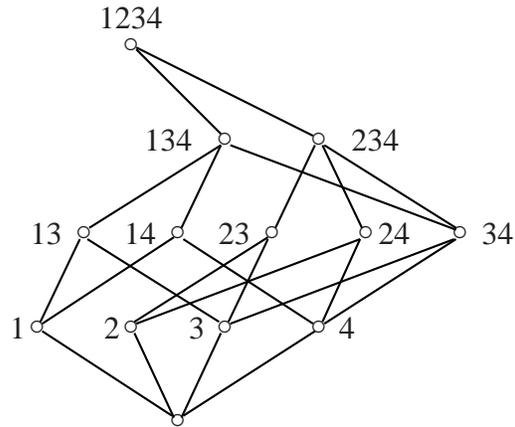
\begin{center}
 \psset{unit=2.5cm}
\begin{pspicture}(-2.5,-0.5)(5.5,2)
\rput(0.75,0){$\circ$} \rput(0,0.5){$\circ$}
\rput(0.5,0.5){$\circ$}
 \rput(1,0.5){$\circ$} \rput(1.5,0.5){$\circ$}
\rput(0.25,1){$\circ$}\rput(0.75,1){$\circ$}
\rput(1.25,1){$\circ$}\rput(1.75,1){$\circ$}
\rput(2.25,1){$\circ$} \rput(1.5,1.5){$\circ$}
\rput(1,1.5){$\circ$}\rput(0.5,2){$\circ$}

\psline(0.72,0.03)(0.03,0.48) \psline(0.73,0.02)(0.52,0.47)
\psline(0.77,0.02)(0.98,0.47)\psline(0.78,0.03)(1.48, 0.48)
\psline(0.02,0.53)(0.23,0.97)\psline(0.03,0.52)(0.73,0.98)
\psline(0.52,0.53)(1.23,0.98)\psline(0.53,0.52)(1.72,0.97)
\psline(0.97,0.53)(0.27,0.97)\psline(1.03,0.52)(2.22,0.97)
\psline(1.02,0.53)(1.23,0.97)\psline(1.47,0.52)(0.78,0.97)
\psline(1.52,0.53)(1.72,0.97)\psline(1.53,0.52)(2.23,0.97)
\psline(0.28,1.03)(0.98,1.47)\psline(0.77,1.03)(0.98,1.47)
\psline(1.27,1.03)(1.48,1.47)\psline(1.73,1.03)(1.52,1.47)
\psline(2.23,1.03)(1.53,1.47)\psline(2.23,1.02)(1.03,1.47)
\psline(0.98,1.53)(0.53,1.98)\psline(1.47,1.52)(0.53,1.99)

\rput(-0.1,0.5){1}\rput(0.4,0.5){2}\rput(0.85,0.5){3}\rput(1.65,0.5){4}
\rput(0.05,1){13}\rput(0.55,1){14}\rput(1.05,1){23}\rput(1.9,1){24}\rput(2.45,1){34}
\rput(0.7,1.5){134}\rput(1.8,1.5){234}\rput(0.5, 2.15){1234}
\end{pspicture}
\end{center}
The $h$-vector of $\LL$ is $(1,4,6,2)$.

\medskip
By Theorem \ref{maintheorem} the $h$-vector of the Alexander dual
a meet-distributive meet-semilattice is just the $f$-vector of a
simplicial complex, and the $h$-vector of a level simplicial
complex coming from a meet-distributive meet-semilattice is just
the $f$-vector of a pure simplicial complex. These facts  lead us
to the following

\begin{Question}{\em  (a) Characterize the $h$-vector of the Alexander
dual of a finite distributive lattice.

(b) Characterize the $h$-vector of the Alexander dual a
meet-distributive lattice.

(c) Find a nice class of level simplicial complexes whose
$h$-vector is  not the $f$-vector of a pure simplicial complex.}
\end{Question}

For example $(1,3,3)$ is the $h$-vector of the Alexander of the
meet-distributive lattice $\BB_{[3]}\setminus\{1,3\}$, but is not
the $h$-vector of the Alexander dual of a distributive lattice.

\end{document}